\documentclass[11pt]{amsart}

\usepackage[a4paper,hmargin=3cm,vmargin=3cm]{geometry}
\usepackage{amsfonts,amssymb,amscd,amstext}
\usepackage{graphicx}
\usepackage[dvips]{epsfig}
\usepackage[ansinew]{inputenc}

\usepackage{fancyhdr}
\usepackage{times}

\usepackage{enumerate}
\usepackage{titlesec}
\usepackage{mathrsfs}
\usepackage{stmaryrd}

\def\R{\mathbb{R}}

\def\D{\mathbb{D}}

\newcommand{\ben}{\begin{enumerate}}
\newcommand{\bit}{\begin{itemize}}
\newcommand{\een}{\end{enumerate}}
\newcommand{\eit}{\end{itemize}}

\newcommand{\ed}{\end{document}}

\def\cU{\mathcal{U}}

\def\cW{\mathcal{W}}

\def\cV{\mathcal{V}}

\def\cM{\mathcal{M}}

\def\cG{\mathcal{G}}
\def\cF{\mathcal{F}}

\def\cO{\mathcal{O}}

\let\8=\infty \let\0=\emptyset

\let\landa=\lambda
\let\alfa=\alpha

\let\parc=\partial

\def\ep{\varepsilon}

\def\landa{\lambda}

\def\flecha{\rightarrow}
\def\esiz{\langle}
\def\esde{\rangle}

\def\cte.{\mathop{\rm cte.}\nolimits}

\def\R{\mathbb{R}}

\def\D{\mathbb{D}}

\def\S{\mathbb{S}}

\titleformat{\section}
{\filcenter\bfseries\large} {\thesection{.}}{0.2cm}{}
\titleformat{\subsection}[runin]
{\bfseries} {\thesubsection{.}}{0.15cm}{}[.]
\titleformat{\subsubsection}[runin]
{\em}{\thesubsubsection{.}}{0.15cm}{}[.]

\usepackage[up,bf]{caption}

\newtheorem{theorem}{Theorem}[section]

\newtheorem{corollary}[theorem]{Corollary}
\newtheorem{definition}[theorem]{Definition}

\newtheorem{example}[theorem]{Example}

\theoremstyle{definition}

\numberwithin{equation}{section}
\numberwithin{figure}{section}

\usepackage{color}

\begin{document}
\fancyhead[LO]{Overdetermined elliptic problems}
\fancyhead[RE]{Pablo Mira}
\fancyhead[RO,LE]{\thepage}

\thispagestyle{empty}

\begin{center}
{\bf \LARGE Overdetermined elliptic problems in topological disks}
\vspace*{5mm}

\hspace{0.2cm} {\Large Pablo Mira}
\end{center}

\footnote[0]{\vspace*{-0.4cm} \emph{Mathematical Subject Classification}: 35J25, 53A10}

\vspace*{7mm}

\begin{quote}
{\small
\noindent {\bf Abstract}\hspace*{0.1cm}
We introduce a method, based on the Poincaré-Hopf index theorem, to classify solutions to overdetermined problems for fully nonlinear elliptic equations in domains diffeomorphic to a closed disk. Applications to some well-known nonlinear elliptic PDEs are provided. Our result can be seen as the analogue of Hopf's uniqueness theorem for constant mean curvature spheres, but for the general analytic context of overdetermined elliptic problems.


\vspace*{0.1cm}

}
\end{quote}


\section{Introduction}
The following famous theorem by Serrin \cite{Se} is widely regarded as the archetypical result on overdetermined problems for elliptic PDEs: \emph{if $u\in C^2(\overline{\Omega})$ solves}
\begin{equation}\label{serrineq}\def\arraystretch{1.6}\left\{\begin{array}{lll} \Delta u+1=0 & \text{ in } & \Omega \\
u=0, \hspace{0.5cm} \displaystyle \frac{\parc u}{\parc \nu} = c  & \text{ on } & \parc \Omega, 
\end{array} \right.
\end{equation} 
\emph{where $\Omega\subset \R^n$ is a bounded smooth open domain, then $\Omega$ is a ball and $u$ is radially symmetric}; here $c$ is constant and $\nu$ is the interior unit normal of $\parc \Omega$. For the proof, Serrin introduced the \emph{method of moving planes}, a boundary version for overdetermined elliptic problems of the geometric Alexandrov theorem \cite{A1}, according to which compact embedded constant mean curvature (CMC) hypersurfaces in $\R^{n+1}$ are round spheres. 

Besides Alexandrov's theorem, the second classical result that models the geometry of compact CMC surfaces is Hopf's theorem \cite{Ho0,Ho}: \emph{compact simply connected CMC surfaces in $\R^3$ are round spheres}. These two theorems are proved by totally different techniques, and complement each other. For instance, while Alexandrov's theorem works for arbitrary dimension and any topological type, Hopf's theorem is specific of dimension two (see Hsiang \cite{Hs} for counterexamples in higher dimension) and needs the surface to be simply connected (Wente tori \cite{W} are counterexamples for the multiply connected case). On the other hand, Hopf's theorem allows arbitrary self-intersections, and its proof provides important information on the local geometry of any CMC surface. Both results have been extremely influential in surface theory.

In this paper we prove what can be seen as a boundary version for overdetermined elliptic problems of Hopf's theorem. Our theorem somehow completes the general parallelism between compact constant mean curvature theory and overdetermined elliptic problems in bounded domains initiated by Serrin in \cite{Se}. For unbounded domains $\Omega$ and semilinear elliptic PDEs, this paralellism has been deeply investigated, see e.g. \cite{HHP,RS,RRS1,RRS2,T}.

The following particular case can be seen as a model situation for our main result (Theorem \ref{mainth}). Consider the overdetermined problem
\begin{equation}\label{prob}\def\arraystretch{1.6}\left\{\begin{array}{lll} F(D^2u,Du)=0 & \text{ in } & \Omega, \\
u=0, \hspace{0.5cm} \displaystyle \frac{\parc u}{\parc \nu} = g(\nu)  & \text{ on } & \parc \Omega, 
\end{array} \right.
\end{equation} 
where $F(D^2u,Du)=0$ is a $C^{1,\alfa}$ fully nonlinear elliptic equation, $\Omega\subset \R^2$ is a bounded $C^2$ domain and $g\in C^1(\S^1)$. We will assume that the pair $(F,g)$ satisfies the following compatibility condition, which ensures that \eqref{prob} has solutions, and that we call \emph{Property $(*)$}:

\begin{quote}
\emph{There is some solution $u^0\in C^2(\Omega_0^*)$ to $F(D^2 u,Du)=0$ whose gradient is an orientation preserving diffeomorphism from $\Omega_0^*$ onto $\R^2$, such that $u^0$ solves \eqref{prob} when restricted to some $C^2$ bounded domain $\overline{\Omega_0} \subset \Omega_0^*$. }
\end{quote}

In this situation, and imposing a boundary regularity condition, we will prove: \emph{if $\Omega$ is simply connected and $u\in C^2(\overline{\Omega})$ solves \eqref{prob}, then up to a translation $u=u^0$ and $\Omega=\Omega_0$}.

There are many well-studied elliptic equations  $F(D^2u,Du)=0$ that satisfy Property $(*)$, and for which a solution to \eqref{prob} was not previously known, not even for planar simply connected domains; examples will be provided in Section \ref{sec:ap}. An analogous theorem holds for elliptic equations of the more general form $F(D^2u,Du,u)=0$, although in that case the corresponding Property $(*)$ is more involved; see Definition \ref{canoni}.

As happens with the classical situation of Alexandrov and Hopf, our theorem complements the previous methods for solving overdetermined elliptic problems (see e.g. \cite{Se,BNST,Wn}). Our method only works for the particular --but fundamental-- case that $\Omega$ is a simply connected planar domain, since it depends on the Poincaré-Hopf index theorem. But on the other hand, it works for fully nonlinear elliptic equations without symmetries (the moving planes method is inapplicable in that situation), and for equations for which no \emph{$P$-function} (in the sense of Weinberger's approach \cite{Wn}) is known. Also, our method does not need to ensure, or impose, that $u>0$ or $u<0$ in $\Omega$, it works for self-intersecting domains (see Subsection \ref{imdi}) and it does not need the solution $u$ to satisfy $F(D^2 u,Du,u)=0$ everywhere on $\Omega$ (see Subsection \ref{umbili}).

We should also observe that the constant overdetermined boundary conditions
\begin{equation}\label{conse}
u=0, \hspace{0.5cm} \frac{\parc u}{\parc \nu} = {\rm constant} \hspace{0.5cm} \text{ on $\parc \Omega$}
\end{equation}
as in \eqref{serrineq} are the natural ones in the case that the PDE $F[u]=0$ is rotationally symmetric (and thus it admits radial solutions), but they lose all interest or prominent role when $F[u]=0$ is not rotationally invariant. In our situation, the existence of the \emph{canonical solution} $u^0$ in Property $(*)$ indicates that the Neumann condition in \eqref{prob} is the natural one; see the discussion in Example \ref{eje1}. Also, note that when $u^0$ is radially symmetric, this natural Neumann condition is the constant one in \eqref{conse}. There are many works on overdetermined elliptic problems with non-constant associated Neumann conditions, see e.g. \cite{BHS,BGHT,CS,ES,O,Se,WC}.

We have organized the paper as follows. In Section \ref{sec:defi} we give some basic definitions about overdetermined elliptic problems and canonical families of solutions to elliptic PDEs of the form $F(D^2u,Du,u)=0$, and we state our main uniqueness result (Theorem \ref{mainth}). 

In Section \ref{sec:main} we prove Theorem \ref{mainth}. The basic idea is to construct, associated to any non-canonical solution to \eqref{prob}, a line field on $\overline{\Omega}$ with isolated singularities of negative index that is tangent to $\parc \Omega$, and to derive from there a contradiction with the Poincaré-Hopf theorem using that $\Omega$ is simply connected. This type of index strategy originates from Cohn-Vossen's proof of the rigidity of ovaloids in $\R^3$, and has been previously used in several ways in surface theory. Our proof also uses ideas from our previous work \cite{GM3} with Gálvez on uniqueness of immersed spheres modeled by elliptic PDEs in three-manifolds, but the proof that we present here is self-contained. 

In Section \ref{sec:ap} we discuss some direct specific applications of our method to well known fully nonlinear elliptic PDEs such as the equations of prescribed mean or Gaussian curvature. We will also indicate how Theorem \ref{mainth} still holds under weaker conditions, like self-intersecting domains, or functions that do not satisfy the PDE at every point. 

It is important to observe that these applications in Section \ref{sec:ap} are just some of the most visible consequences of the general index method that we present here, which seems suitable to study in great generality elliptic problems in simply connected planar domains. For example, this index method can be used to classify solutions to fully nonlinear anisotropic problems, in the spirit of \cite{CS,WC}; we will explain this briefly in Section \ref{sec:ap}. The method is also specially suitable to study nonlinear elliptic equations that might depend on the independent variables $(x,y)$, but which are invariant with respect to some Lie group structure in $\R^3$ different from its canonical abelian structure. However, for the sake of brevity, these lines of inquiry will not be pursued here.

The author is grateful to J.A. Gálvez, A. Enciso and D. Peralta-Salas for helpful discussions.
\section{Statement of the theorem}\label{sec:defi}

Throughout this paper, we will let $F\in C^{1,\alfa} (\cU)$, where $\cU=\R^3\times \cV\subset \R^6$ with $\cV$ a convex open set, and so that $F=F(z,p,q,r,s,t)$ satisfies on $\cU$ the ellipticity condition $4 F_r F_t - F_s^2>0$. We will denote by $F[u]=0$ the fully nonlinear elliptic PDE in two variables
  \begin{equation}\label{fulpde}
 F(u,u_x,u_y,u_{xx},u_{xy},u_{yy})=0.
 \end{equation}
\begin{definition}\label{canoni}
A \emph{canonical family of solutions} to \eqref{fulpde}  is a family $\cF=\{u_t : t\in \R\}$ of $C^2$ solutions to \eqref{fulpde} with the following properties:
 \begin{enumerate}
 \item
$D u_t : \Omega_t \flecha \R^2$ is an orientation preserving diffeomorphism for every $u_t:\Omega_t\flecha \R$ in $\cF$.
 \item
The family $\cF$ is $C^1$ with respect to the parameter $t$.
 \item
For every $(z,p,q)\in \R^3$ there exist a unique $t\in \R$ such that $u_t(x,y)=z$ and $Du_t(x,y)=(p,q)$ for some $(x,y)\in \Omega_t$.
 \end{enumerate}
\end{definition}
When \eqref{fulpde} is of the form $F(D^2u,Du)=0$, if there exists a solution $u^0\in C^2 (\Omega_0^*)$ to \eqref{fulpde} satisfying that $Du^0 : \Omega_0^*\flecha \R^2$ is an orientation preserving diffeomorphism,
%
then the family $\{u_t := u^0+t : t\in \R\}$ is a canonical family of solutions to \eqref{fulpde}. For example, the family 
$u_t(x,y):= x^2+y^2+t$ is a canonical family of solutions for the Monge-Ampère equation ${\rm det}(D^2 u) =4$ in $\R^2$.

Let now $u^0=u_t\in \cF$ be a canonical solution for which $\gamma:=(u^0)^{-1}(0)$ is a regular curve contained in $\Omega=\Omega_t$. Then, $\gamma$ is a strictly convex regular curve, and if $\nu$ denotes the inner unit normal vector along $\gamma$ associated to the convex domain $\Omega\subset \R^2$ with $\parc \Omega=\gamma$, then there exists $g\in C^1(\S^1)$ such that  $u^0$ solves the overdetermined problem 
\begin{equation}\label{natu0}
\def\arraystretch{1.6}\left\{\begin{array}{lll} F[u]=0 & \text{ in } & \Omega, \\
u=0, \hspace{0.5cm} \displaystyle \frac{\parc u}{\parc \nu} = g(\nu)  & \text{ on } & \parc \Omega, 
\end{array} \right.
\end{equation} 

Note that if $u^0$ is radially symmetric with respect to some point $p_0\in \R^2$, the associated function $g$ is constant, and we recover the classical overdetermined conditions \eqref{conse}.

\begin{definition}
We will call problem \eqref{natu0} the \emph{natural overdetermined problem} associated to the canonical solution $u^0=u_t\in \cF$  to $F[u]=0$. Note that each family $\cF$ gives rise to a one-parameter family of problems \eqref{natu0}, one for each value $t\in \R$ where $u_t^{-1}(0)$ is a regular $C^2$ curve contained in $\Omega_t$.
\end{definition}


\begin{example}\label{eje1}
\emph{The following discussion might be illustrative of why we call \eqref{natu0} a natural overdetermined problem. After a change $(x,y)\mapsto (x,2y)$, the linear equation $\Delta u +1=0$ in dimension two changes to }
 \begin{equation}\label{cha}
 u_{xx}+4u_{yy}+1=0,
 \end{equation}
\emph{which is no longer rotationally invariant. The radial solutions $u(x,y)=a-\frac{x^2+y^2}{4}$ to $\Delta u +1 =0$ are transformed into $u_a(x,y)= a-\frac{x^2}{4} - \frac{y^2}{16}$, which are the simplest solutions to \eqref{cha}. These solutions intersect the $z=0$ plane when $a>0$ along an ellipse, but the intersection angle is not constant anymore, i.e. $u_a$ does not satisfy the boundary conditions \eqref{conse}. Moreover, these overdetermined constant conditions \eqref{conse} for \eqref{cha} lose all interesting meaning, and in general will not support a solution. Still, by the convexity of $u_a(x,y)$ and of the intersection ellipse $\Omega_a$, there is a certain function $g_a:\S^1\flecha \R$, which can be explicitly calculated, such that $\frac{\parc u_a}{\parc \nu} = g_a(\nu)$ along $\parc \Omega_a := u_a^{-1}(0)$. Thus, we can formulate the \emph{natural} overdetermined problem for \eqref{cha} as }
\begin{equation}\label{serrino}\def\arraystretch{1.8}\left\{\begin{array}{lll} u_{xx}+4u_{yy}+1=0 & \text{{\emph in} } & \Omega, \\
u=0, \hspace{0.5cm} \displaystyle \frac{\parc u}{\parc \nu} = g_a(\nu)  & \text{\emph on } & \parc \Omega, 
\end{array} \right.
\end{equation} 
\emph{which will have $u_a$ as the solution for which uniqueness is aimed.}
\end{example}

Our main theorem is a general uniqueness result for natural overdetermined problems in the case that $\overline{\Omega}$ is diffeomorphic to a closed disk. Specifically, we will consider problem \eqref{natu0}, where
 \begin{enumerate}
 \item
$F[u]=0$ is the fully nonlinear equation \eqref{fulpde}.
 \item
There is a canonical family of solutions $\cF=\{u_t : t\in \R\}$ to \eqref{fulpde}.
 \item
$g\in C^1(\S^1)$ is a \emph{natural Neumann condition}, i.e. it is given by $\frac{\parc u^0}{\parc \nu} =g(\nu)$ for some $u^0\in \cF$ that intersects the $z=0$ plane along a closed convex curve $\gamma$; here $\nu$ is the inner unit normal to $\gamma$.
 \end{enumerate}
In these conditions, we have:

\begin{theorem}\label{mainth}
Let $u\in C^2(\Omega_{\ep})$ denote a solution to $F[u]=0$ that solves \eqref{natu0}
when restricted to some compact $C^2$ simply connected domain $\overline{\Omega}\subset \Omega_{\ep}$. Then:
\begin{enumerate}
 \item
$\Omega$ is a translation of the domain $\Omega_0$ bounded by $\gamma$, i.e. $\Omega$ is of \emph{canonical shape}.
\item
$u=u^0\circ T$ for some translation $T$ of $\R^2$.
\end{enumerate}
\end{theorem}

We should observe that we are making a strong boundary regularity assumption in Theorem \ref{mainth}, by requiring that $u$ solves $F[u]=0$ on an open domain $\Omega_{\ep}$ containing $\overline{\Omega}$. In the real analytic case, this assumption is not necessary. In the general case, Theorem \ref{mainth} likely holds with far less restrictive boundary regularity assumptions, but it is not the purpose of this short note to discuss this aspect.



\section{Proof of Theorem \ref{mainth}}\label{sec:main}
Note that, by interior regularity \cite{N2}, any $C^2$ solution to \eqref{fulpde} is of class $C^{3,\beta}$, $0<\beta<1$.
Consider the family 
 \begin{equation}\label{fag}
\cG:=\{u_t (x+a,y+b) : u_t\in \cF, t\in \R, (a,b)\in \R^2\}.
\end{equation}
By the last two properties in Definition \ref{canoni} it is clear that $\cG$ depends $C^1$ smoothly on $(a,b,t)$, and that for every $\alfa=(x,y,z,p,q)\in \R^5$ there exists a unique $u^{\alfa} \in \cG$ such that $$\alfa=(x,y,u^{\alfa}(x,y),Du^{\alfa}(x,y))$$ for some $(x,y)\in {\rm dom}(u^{\alfa})$. Thus, we can write $\cG=\{u^{\alfa} : \alfa\in \R^5\}$, and this family is $C^1$ with respect to $\alfa$. In particular,
 \begin{equation}\label{gamap}
 \Gamma(\alfa) := (D^2 u^{\alfa})_{(x,y)}, \hspace{1cm} \alfa=(x,y,z,p,q),
 \end{equation}
defines a $C^1$ map from $\R^5$ into the space of positive definite symmetric bilinear forms in $\R^2$.

Let now $u\in C^2(\Omega_{\ep})$ be as in the statement, and consider the map $$\Lambda: (x,y)\in \Omega_{\ep} \mapsto \Gamma((x,y,u(x,y),Du(x,y)).$$ Clearly, $\Lambda$ defines a $C^1$ Riemannian metric on $\Omega_{\ep}$. This allows to define a $C^1$ tensor $S$ on $\Omega_{\ep}$ given by 
 \begin{equation}\label{def:s}
 \Lambda (S(X),Y) =D^2 u (X,Y)
 \end{equation}
for every pair $X,Y$ of tangent vector fields on $\Omega_{\ep}$; the existence of such $S$ follows since $D^2 u$ is symmetric and $\Lambda$ is positive definite. As $S$ is diagonalizable, at every $p\in \Omega_{\ep}$ where $S_p$ is not proportional to the identity there exist exactly two eigenlines for $S_p$, both of them orthogonal with respect to the Riemannian metric $\Lambda$ at $p$. Moreover, it is trivial to observe from \eqref{def:s} that, at an arbitrary point $p\in \Omega_{\ep}$, $S_{p}(w)=\landa w$ for some $w\neq 0$ if and only if \begin{equation}\label{carbo}
\landa (D^2 u^{\xi})_{p} (w,Y)= (D^2 u)_p (w,Y)\end{equation} for every $Y\in \R^2$, where $\xi=(p,u(p),Du(p))$.

Assume now that $p\in \Omega_{\ep}$ is a point where $S_p=\landa {\rm Id}$ for some $\landa \in \R$. By \eqref{carbo}, this means that 
 \begin{equation}\label{compla}
\landa (D^2 u^{\xi})_p = (D^2 u)_p ,\end{equation}
where $\xi=(p,u(p),D(p))$. But as $u^{\xi}$ and $u$ are both solutions to $F[u]=0$, the ellipticity of $F$ implies that \eqref{compla} cannot happen unless $\landa=1$, i.e. unless $S_p={\rm Id}$ and $D^2(u-u^{\xi})_p=0$.

As a result, if $\cO:=\{p\in \Omega_{\ep} : S_p={\rm Id}\}$, we can define on $\Omega_{\ep}\setminus \cO$ a pair of $C^1$ line fields $Z_1,Z_2$, both of them orthogonal with respect to the Riemannian metric $\Lambda$, and given by the eigenlines of $S$. 
Our next aim is to analyze the behavior of these line fields $Z_1,Z_2$ around their singularities, i.e. around points in $\cO$.

Let $p=(x_0,y_0)\in \cO$, and denote again $\xi=(p,u(p),D(p))$. Then, $F[u]=F[u^{\xi}]=0$. Denote $\phi=u-u^{\xi}$. By ellipticity of $F$, there exists an elliptic linear homogeneous operator of second order $L$ with $C^{\alfa}$ coefficients such that $L[\phi]=0$. Also, observe that $\phi(p)$, $D\phi(p)$, $(D^2\phi)_p$ all vanish. Hence, by Bers' theorem \cite{Be}, either $\phi=0$ around $p$, or there exists a homogeneous polynomial $h$ of degree $n\geq 3$ such that 
 \begin{equation}\label{asi1}
 \phi(q)= h(q-p) + o (|q-p|)^n
 \end{equation}
for $q$ sufficiently close to $p$, and so that $h$ is harmonic with respect to some coordinates $(x',y')$ obtained from $(x,y)$ after an affine change of variables. Observe that if $\phi$ is identically zero, then $p\in {\rm int}(\cO)$. 

Consider next some $p \not\in {\rm int}(\cO)$. So, \eqref{asi1} holds. Define a $C^1$ symmetric bilinear form $\sigma$ on $\Omega_{\ep}$ by $\sigma=  D^2 u - \Lambda$. Note that for any $q\in \Omega_{\ep}$ we have 
 \begin{equation}\label{neq1} \sigma(q)= (D^2 \phi)_q + (D^2 u^{\xi})_q - \Lambda (q).
  \end{equation}
Using \eqref{asi1} and the definition of $\Lambda$, we obtain from \eqref{neq1}
 \begin{equation}
 \sigma(q)=(D^2 h)_{q-p} + \Gamma(q,u^{\xi} (q),Du^{\xi} (q)) - \Gamma(q,u(q),Du(q)) + o (|q-p|)^{n-2},
 \end{equation}
which by the mean value theorem applied to $\Gamma$, reduces by \eqref{asi1} to
 \begin{equation}\label{asi2}
 \sigma(q)=(D^2 h)_{q-p} + o (|q-p|)^{n-2}.
 \end{equation}
In particular, since $h$ is a harmonic homogeneous polynomial of degree $n\geq 3$ for the coordinates $(x',y')$, it follows that ${\rm det} (\sigma(q))<0$ for any $q\in \Omega_{\ep}$ in a sufficiently small punctured neighborhood of the point $p\in \cO$. This implies, in particular, that $p$ is isolated as an element of $\cO$, and so that $p$ is an isolated singularity of the two line fields $Z_1,Z_2$.

This proves that $\cO\subset \Omega_{\ep}$ is composed only by interior points (they correspond to the case where $\phi$ vanishes identically around the point) or by isolated points. A simple topological argument ensures then that either $\cO=\Omega_{\ep}$, or $\cO\cap \overline{\Omega}$ is a finite (posibly empty) set. Moreover, in the case $\cO=\Omega_{\ep}$ we clearly have from the previous argument that $u=u^{\xi}$ for any $\xi\in \R^5$ of the form $\xi=(p,u(p),Du(p))$ for some $p\in \Omega_{\ep}$.

Next, we will assume that $\cO\neq \Omega_{\ep}$ and obtain a contradiction.

To start, we will show that one of $Z_1,Z_2$ is tangent to $\parc \Omega$. First, let us prove that the tangent line to $\parc \Omega$ is an eigenline of $S$ for every $p\in \parc \Omega$. By previous arguments, this property is equivalent to showing that for each $p\in \parc \Omega$ there exists some $\landa=\landa(p)$ so that \eqref{carbo} holds for every $Y \in \R^2$, where $w$ is the positively oriented unit tangent vector to $\parc \Omega$ at $p$.

Let $\gamma(s)$ be a parametrization of $\parc \Omega$ with $\gamma(0)=p$, $\gamma'(0)=w$. The first and second derivatives of equation $u(\gamma(s))=0$ evaluated at $s=0$ give, respectively,
 \begin{equation}\label{be1}
 \esiz w, Du(p)\esde =0 \hspace{0.5cm} \text{and} \hspace{0.5cm} \esiz \gamma''(0), Du(p)\esde + (D^2 u)_p (w,w)=0.
 \end{equation}
Denote $\nu_0:=\nu(0)$, where $\nu(s)$ is the inner unit normal along $\gamma(s)$.
Since the Neumann condition gives $\esiz \nu_0,Du(p)\esde = g(\nu_0)$, we obtain from \eqref{be1}
 \begin{equation}\label{be3}
 (D^2 u)_p (w,w)= - \kappa (p) g(\nu_0),
 \end{equation}
where $\kappa(p)$ denotes the curvature of $\parc \Omega$ at $p$. 

Let now consider the element $u^{0}$ of the canonical family $\cF$ that solves \eqref{natu0}. After a translation $T:(x,y)\mapsto (x+a,y+b)$, $u^0\circ T$ coincides with the element $u^{\xi}\in \cG$ given by the choice $\xi=(p,0,Du(p))$. Once here, the same computation applied to $u^{\xi}$ shows that 
 \begin{equation}\label{be4}
  (D^2 u^{\xi})_p (w,w)= - \kappa^{\xi} (p) g(\nu_0),
 \end{equation}
 where $\kappa^{\xi}(p)$ is the curvature at $p$ of the boundary $\parc\Omega^{\xi}:=(u^{\xi})^{-1}(0)\subset \R^2$, which by strict convexity of $u^{\xi}$ is always positive. Thus, we see that \eqref{carbo} holds at $p$ for $Y=w$ and 
 \begin{equation}\label{be6}
 \landa= \frac{\kappa(p)}{\kappa^{\xi} (p)}.
 \end{equation}

Next, if we differentiate the Neumann condition $\esiz \nu(s),Du(\gamma(s))\esde =g(\nu(s))$ and evaluate it at $s=0$, a similar computation shows that 
 \begin{equation}\label{be5}
 (D^2 u)_p (w,\nu_0) = - \kappa(p) (dg)_{\nu_0} (w), \hspace{1cm} (D^2 u^{\xi})_p (w,\nu_0) = - \kappa^{\xi}(p) (dg)_{\nu_0} (w).
 \end{equation}
Putting together \eqref{be6}, \eqref{be5} we see that \eqref{carbo} holds at $p$ for $Y=\nu_0$ and $\landa$ given by \eqref{be6}. So, by linearity, \eqref{carbo} holds for this special value of $\landa$ and for all $Y\in \R^2$. Hence, the tangent line to $\parc \Omega$ is an eigenline of $S$ at every $p\in \parc \Omega$. But now let us recall that:
 \begin{enumerate}
 \item
The eigenlines $Z_1,Z_2$ are well defined, unique and continuous at every $p\in \parc \Omega - \mathcal{O}$, and $\mathcal{O}\cap \parc \Omega$ is a finite set.
 \item
At every $p\in \parc\Omega -\mathcal{O}$, $Z_1$ and $Z_2$ are orthogonal with respect to the metric $\Lambda$.
 \end{enumerate}
Since $\parc \Omega$ is differentiable, this implies that one of the line fields $Z_1,Z_2$ is everywhere tangent to $\parc \Omega$, as we wished to show.

Let now $p$ be any point in $\cO$; by the previous arguments, $p$ is an isolated point of $\cO$ and \eqref{asi2} holds. The fact that $\sigma$ is indefinite in a punctured neighborhood of $p$ allows to consider for each $q$ in that punctured neighborhood the \emph{null lines} $(U,V)$ of $\sigma_q:=\sigma(q)$, given by $\sigma_q(U,U)=\sigma_q (V,V)=0$. Similarly, we can consider for any such $q$ the null lines $(U^{h},V^{h})$ of $(D^2 h)_q$. In both cases, these null lines define two continuous line fields around $p$ with an isolated singularity at $p$. Since $h$ is a homogeneous harmonic polynomial of degree $n\geq 3$ in the $(x',y')$ coordinates, a standard computation shows that the index around $p$ of  the null directions of the Hessian of $h$ with respect to $(x',y')$, which are given by ${\rm Re} (h_{\zeta \zeta} d\zeta^2)=0$ where $\zeta=x'+iy'$, is 
equal to $-(n-2)/2<0$. As $(x',y')$ differ from $(x,y)$ by an affine transformation of $\R^2$, we conclude that the index of both $U^{h}$ and $V^{h}$ around $p$ is negative. Finally, by \eqref{asi2}, the index of the line fields $U$ and $V$ around $p$ is also negative.

Let now $\{e_1,e_2\}$ be a basis of eigenvectors of $S_q$. As they are orthogonal with respect to $\Lambda$, a simple computation shows that if $\sigma_q(w_1,w_1)=0$ for $w_1=x_1e_1+x_2e_2$, then $\sigma_q (w_2,w_2)=0$ for $w_2:=x_1 e_2 - x_2 e_2$. That is, the eigenlines $Z_1,Z_2$ of $S$ bisect (with respect to $\Lambda$) the null line fields $U,V$ at every point around $p$. Thus, $Z_1,Z_2$ also have negative index around $p$. 

In this way, we have created a continuous line field $Z$ on $\Omega_{\ep}$ with only isolated singularities, all of them of negative index, and which is tangent to $\parc \Omega$. Thus, the \emph{boundary index} of $Z$ as a line field on $\overline{\Omega}$ at these isolated singularities that lie in $\parc \Omega$ is one half of their index as singularities of $Z$ in $\Omega_{\ep}$. Since $Z$ is tangent along $\parc \Omega$, after topologically identifying $\overline{\Omega}$ with the closed upper hemisphere $S^+:=\overline{\S_+^2}$, we can extend $Z|_{S^+}$ to a continuous line field on $\S^2$, by reflecting $Z$ across $\parc S^+$ in a symmetric way. In this way, we obtain a continuous line field in $\S^2$ with a finite number of singularities, all of them of negative index. Since the Poincaré-Hopf theorem implies that the sum of the indices of a continuous line field with isolated singularities in $\S^2$ is equal to $2$ and this extended line field does not satisfy this property, we reach a contradiction.

Consequently, $\cO=\Omega_{\ep}$, and the function $\phi=u-u^{\xi}$ vanishes identically on $\Omega_{\ep}$ for any $\xi\in \R^5$ of the form $\xi=(p,u(p),Du(p))$ for some $p\in \Omega_{\ep}$. But now, let us recall that the element $u^0\in \cF$ that solves \eqref{natu0} satisfies $u^0\circ T =u^{\xi}$ for some translation $T$ of $\R^2$ and some $\xi$ of the previous form. The statement of the theorem follows then trivially.
  
 \section{Discussion of the result: applications and extensions}\label{sec:ap}
 
 \subsection{Application to elliptic equations invariant by translations}
Theorem \ref{mainth} is specially useful in the case that \eqref{fulpde} does not depend on $u$, i.e. when it is an elliptic equation of the form $F(D^2u,D u)=0$. For this equation, as explained in Section \ref{sec:defi}, if there exists one solution $u^0\in C^2(\Omega^*_0)$ such that: 
 \begin{enumerate}
 \item
$Du^0:\Omega^*_0\flecha \R^2$ is an orientation preserving diffeomorphism, and
 \item
$\gamma:=(u^0)^{-1} (0)$ is a regular curve in the $z=0$ plane,
\end{enumerate}
then we can define the \emph{natural overdetermined problem} for $F(Du,D^2u)=0$ given by the boundary conditions
 \begin{equation}\label{overda}
u=0, \hspace{0.5cm} \frac{\parc u}{\parc \nu}= g(\nu) \hspace{0.5cm} \text{ on $\parc \Omega$,}\end{equation} where $g\in C^1(\S^1)$ is the function determined by $\frac{\parc u^0}{\parc \nu} = g(\nu)$ along $\gamma$. So, in these conditions, Theorem \ref{mainth} directly implies:

\begin{corollary}\label{noved}
Let $u\in C^2(\Omega_{\ep})$ be a solution to $F(Du,D^2u)=0$ with boundary data \eqref{overda} on some compact simply connected $C^2$ subset $\overline{\Omega}\subset \Omega_{\ep}$. Then (up to a translation) $u=u^0$ and $\Omega=\Omega_0$, where $\Omega_0\subset \R^2$ is the domain bounded by $\gamma=(u^0)^{-1}(0)$.
\end{corollary}
 
There are many well-known nonlinear elliptic equations of the form $F(Du,D^2 u)=0$ for which one can ensure the existence of a canonical solution $u^0$ as above, and for every such PDE we can solve its associated natural overdetermined problem on compact simply connected domains via Corollary \ref{noved}. We discuss next some examples.

\vspace{0.2cm}
 {\it {\bf 1.} The Monge-Ampère equation of Minkowski's problem}, given by
\begin{equation}\label{mami}{\rm det} (D^2 u) = W (Du) (1+|Du|^2)^2,\end{equation}
where $W\in C^2(\R^2)$ is given by 
 \begin{equation}\label{wf}
f(x_1,x_2,x_3)=W\left(\frac{x_1}{x_3},\frac{x_2}{x_3}\right) \hspace{1cm} \forall (x_1,x_2,x_3)\in \S_+^2
\end{equation}
for a function $f>0$ in $\S^2$ such that $\int_{\S^2} x/f(x) \, dx =0$; by classical works (cf. \cite{N,Po}, see also \cite{Le}), this integral condition ensures the existence of the canonical solution $u^0$ to \eqref{mami}.

\vspace{0.2cm}

{\bf 2.} Similary, we can consider the \emph{prescribed mean curvature equation}
\begin{equation}\label{pme}
{\rm div}\left(\displaystyle\frac{Du}{\sqrt{1+|Du|^2}}\right) = 2 W (Du),\end{equation}
where $W$ is given by \eqref{wf} for $f\in C^2(\S^2)$ with $f(-x)=f(x)>0$. In these conditions, seminal work by B. Guan and P. Guan \cite{GG} on existence of prescribed mean curvature ovaloids in $\R^{n+1}$ ensures that a canonical solution to \eqref{pme} exists. A general uniqueness theorem for Guan-Guan spheres in $\R^3$ can be found in \cite{GM}. 

\vspace{0.2cm}

There are other curvature equations similar to \eqref{mami} or \eqref{pme} that also admit this discussion, e.g. the elliptic equation corresponding to Christoffel's problem. Furthermore, as mentioned in the introduction, some functions based on the \emph{Wulff shape} associated to a smooth norm $H=H(\xi_1,\xi_2)$ in $\R^2-\{0\}$ act as the canonical example $u^0$ of fully nonlinear anisotropic elliptic equations of the form $\Phi(Q_{H} [u],R_{H} [u])=0$, where $Q_H[u],R_H[u]$ denote, respectively, the trace (\emph{anisotropic Laplacian}) and determinant of the anisotropic functional $$\cM_H [u] := D_{\xi}^2 V(Du)D^2 u, \hspace{0.5cm} \text{ where } \hspace{0.5cm} V(\xi_1,\xi_2):=\frac{1}{2} H(\xi_1,\xi_2)^2,$$ see e.g. \cite[pg. 870]{CS}. Corollary \ref{noved} lets us solve then the natural overdetermined problem associated to these fully nonlinear anisotropic equations on compact simply connected planar domains $\Omega\subset \R^2$, proving in particular that any such solution domain $\Omega$ is \emph{of Wulff shape}. For the anisotropic version $Q_H[u]=-1$ of Serrin's theorem, see \cite{CS,WC}.




 \subsection{Extension to immersed disks}\label{imdi}
Problem \eqref{natu0} is also meaningful for immersed domains. For example, let $\psi:\overline{\D}\flecha \R^2$ be a $C^2$ immersion of the closed disk into $\R^2$, and let $\overline{\Omega} =\psi(\overline{\D})$ denote the corresponding immersed disk. We say that $v\in C^2(\overline{\D})$ is a solution to $F[u]=0$ on the immersed disk $\overline{\Omega}$ if for every $p\in \overline{\D}$ the function $u$ given locally around $\psi(p)$ by $v=u\circ \psi$ is a solution to $F[u]=0$. If this PDE is of the form \eqref{fulpde} and admits a canonical family of solutions $\cF=\{u_t : t\in \R\}$, then the boundary conditions \eqref{overda} clearly make sense for this immersed case, and the same proof of Theorem \ref{mainth} works. In this way, we have:

\begin{corollary}\label{ime}
Theorem \ref{mainth} also holds when $\Omega$ is an immersed disk in $\R^2$.
\end{corollary}

This extension is somehow related to a classical work by Nitsche \cite{Ni}, who proved that any constant mean curvature disk immersed in $\R^3$ that intersects a sphere along its boundary at a constant angle is a flat disk or a spherical cap. For that, he used that the Hopf differential of a CMC surface in $\R^3$ is a holomorphic quadratic differential that, in suitable conformal parameters, is real along curvature lines of the surface.

 \subsection{The equation at $\cF$-umbilics}\label{umbili}
 Let $\cF=\{u_t : t\in \R\}$ be a family of canonical solutions to \eqref{fulpde}, and consider the family $\cG$ defined in \eqref{fag}. Define the set
  \begin{equation}\label{defw}
  \cW_{\cF}=\{(x,y,u^{\alfa}(x,y),Du^{\alfa}(x,y),\landa D^2 u^{\alfa}(x,y)) : u^{\alfa}\in \cG, (x,y)\in {\rm dom}(u^{\alfa}), \landa \in \R\},
  \end{equation}
 which is a $6$-dimensional subset of $\R^8$. For a $C^2$ function $u(x,y)$, we will say for short that $u(p_0)\in \cW_{\cF}$ if $(p_0,u(p_0),Du(p_0),D^2u(p_0))\in \cW_{\cF}$. Thus, $u(p_0)\in \cW_{\cF}$ if and only if the Hessian of $u$ at $p_0$ is proportional to the Hessian of $u^{\xi}$ at $p_0$, where $u^{\xi}$ is the unique element of $\cG$ such that $(u(p_0),Du(p_0))=(u^{\xi}(p_0),Du^{\xi}(p_0))$. The existence and uniqueness of this $u^{\xi}$ was discussed in Section \ref{sec:main}. Motivated by classical differential geometry, we introduce the following notion:
 
 \begin{definition}
 We say that a $C^2$ function $u(x,y)$ has an $\cF$-umbilic at $p_0$ if $u(p_0)\in \cW_{\cF}$.
 \end{definition}
With this definition in mind, and going through the proof of Theorem \ref{mainth}, it can be checked that the fact that $u\in C^2(\Omega)$ satisfies \eqref{fulpde} is only used around the $\cF$-umbilics of $u$. Thus, the following more general statement holds.

\begin{corollary}\label{umbi}
Theorem \ref{mainth} holds when $u\in C^2(\Omega_{\ep})$ satisfies the elliptic equation \eqref{fulpde} in some open set $A\subset \Omega_{\ep}$ that contains all $\cF$-umbilics of $u$ (but maybe not globally in $\Omega_{\ep}$).
\end{corollary}
We point out that Corollary \ref{umbi} cannot be deduced by moving planes methods, even in the most symmetric cases like Serrin's $\Delta u +1=0$.
 
\def\refname{References}

\vskip 0.4cm

\noindent Pablo Mira

\noindent Departamento de Matemática Aplicada y Estadística,\\ Universidad Politécnica de Cartagena (Spain).

\noindent  e-mail: {\tt pablo.mira@upct.es}

\vskip 0.4cm

\noindent Research partially supported by MICINN-FEDER,
Grant No.  MTM2013-43970-P, and Programa de Apoyo a la Investigacion,
Fundación Séneca-Agencia de Ciencia y Tecnologia
Region de Murcia, reference 19461/PI/14.

\end{document}